\documentclass[12pt]{article}   
\usepackage{graphicx,psfrag}
\usepackage{color}
\usepackage{array}
\usepackage{colortbl}
\usepackage{tabularx}
\newtheorem{theorem}{Theorem} [section]
\usepackage{graphicx}   
\newtheorem{corollary}[theorem]{Corollary}

\usepackage{amsmath}
\usepackage{amsfonts}
\newfont{\bb}{msbm10}

\def\:{\! :\!}

\usepackage{tikz}

\begin{document}

\title{Cyclic Matching Sequencibility of Graphs}

 \author{Richard A. Brualdi, Kathleen P. Kiernan, Seth A. Meyer\\
 Department of Mathematics\\
 University of Wisconsin\\
 Madison, WI 53706\\
 {\tt \{brualdi,kiernan,smeyer\}@math.wisc.edu}
 \and
 Michael W. Schroeder\\
 Department of Mathematics\\ 
 Marshall University\\
 Huntington, WV\\
 {\tt schroederm@marshall.edu}
 }

\maketitle

 \begin{abstract} 
 We define the cyclic matching sequencibility of a graph to be the largest integer $d$ such that there exists a cyclic ordering of its edges so that every $d$ consecutive edges in the cyclic ordering form a matching. We show that the cyclic matching sequencibility of $K_{2m}$ and $K_{2m+1}$ equal  $m-1$.


\end{abstract}

\section{Introduction}

Let $G\subseteq K_n$ be a graph of order $n$ (without loops or multiple edges) with $m$ edges.   Alspach \cite{Al} defined the following
notions. The {\it matching number} of a linear ordering 
$[1],[2],\ldots,[m]$ of the edges of $G$ (we use $[i]$ to denote the $i$th edge in an ordering of the edges of $G$) is the largest number $d$ such that every $d$ consecutive edges in the ordering form a matching of $G$.
The {\it matching sequencibility} of $G$, denoted ${\rm ms}(G)$, is the maximum matching number of a linear ordering of the edges of $G$. Clearly, the matching sequencibility of $G$ is bounded above
by the largest number of edges in a matching of $G$ and this in turn is bounded above by
$\lfloor \frac{n-1}{2}\rfloor$ if $n$ is odd. If $n$ is even and $G$ itself is a perfect matching, then ${\rm ms}(G)=\frac{n}{2}$. If $n$ is even and $G$ is not a matching, then ${\rm ms}(G)$ cannot equal $ \frac{n}{2}$ because the  edge $[1]$ and the edge $\left[\frac{n}{2}+1\right]$  would have to be identical.  Thus, if $n$ is even, provided $G$ is not a matching, ${\rm ms}(G)\le \lfloor\frac{n-1}{2}\rfloor$.  

The matching sequencibility of $K_2$ is clearly 1. Using the Walecki decomposition of the complete graph $K_n$ with odd $n\ge 3$ 
into Hamilton cycles
and the decomposition of $K_n$ with even $n\ge 4$  into Hamilton paths, Alspach  \cite{Al} showed how to order the edges of $K_n$ ($n\ge 3$) to get a matching number equal to $\lfloor \frac{n-1}{2}\rfloor$. Thus
\begin{equation}\label{eq:alspach}
{\rm ms}(K_n)=\left\lfloor \frac{n-1}{2}\right\rfloor.\end{equation}

Let $G\subseteq K_{p,q}$ be a bipartite graph  with $m$ edges., and let  $A$ be a  $p\times
q$ biadjacency matrix of $G$. Then a ordering of the $m$ edges of
$G$ corresponds to a  bijective replacement of the $m$\ 1s of $A$ with
the integers $1,2,\ldots,m$ resulting in a matrix $\widehat{A}$.  The  
matching sequencibility of this ordering is the largest integer $k$ such
that every  set of $k$ consecutive integers lie in different
rows and columns.

For a complete bipartite graph $K_{m,n}$ with $m\le n$, it is not difficult to show that 
\[{\rm ms}(K_{m,n}) =\left\{\begin{array}{cl}
n-1&\mbox{ if $m=n$}\\
m&\mbox{ if $m<n$.}\end{array}\right.\]
The matching sequencibility is certainly bounded above by these numbers.
Using the biadjacency matrix $A$, it is straightforward to order the edges
(the 1s of $A$) to show that  these upper bounds can be attained. In this representation, a matching of size $k$ corresponds to a set of $k$\ 1s  of $A$ no two from the same row and column. We illustrate this for   $K_{4,4}$ and $K_{4,6}$  whose biadjacency matrices are  the $4\times 4$ and $4\times 6$ matrices $J_{4,4}$ and $J_{4,6}$ of all 1s, respectively. 
\[(K_{4,4}): \widehat{A}=\left[\begin{array}{cccc}
1&5&9&13\\
14&2&6&10\\
11&15&3&7\\
8&12&16&4\end{array}\right]\mbox{ and }
(K_{4,6}): \widehat{A}=\left[\begin{array}{cccccc}
1&13&9&21&5&17\\
18&2&14&10&22&6\\
7&19&3&15&11&23\\
24&8&20&4&16&12\end{array}\right].\]

Henceforth, to avoid trivialities, we assume that $n\ge 3$.
We define a variation of the matching sequencibility of $G$  by considering cyclic orderings of the edges of $G$.
The {\it  cyclic matching number} of a cyclic ordering 
$[1],[2],\ldots,[m],[1]$ of the edges of $G$ is the largest number $d$ such that every $d$ consecutive edges in the cyclic  ordering form a matching of $G$.
The {\it cyclic matching sequencibility} of $G$, denoted ${\rm cms}(G)$, is the maximum cyclic matching number of a cyclic  ordering of the edges of $G$. We have
\begin{equation}\label{eq:ineq}
{\rm cms}(G)\le {\rm ms}(G) \le \left\lfloor\frac{n-1}{2}\right\rfloor.\end{equation}

Clearly, ${\rm cms}(K_3)=1$.
The main result of this note is the following theorem.

\begin{theorem}\label{th:main}
For $n\ge 4$,
\[{\rm cms}(K_n)=
\left\{\begin{array}{cl}
\lfloor \frac{n-1}{2}\rfloor&\mbox{ if $n$ is even}\\
\lfloor \frac{n-3}{2}\rfloor&\mbox{ if $n$ is odd.}\end{array}\right.\]  
\end{theorem}

Combining Theorem \ref{th:main} and equation  (\ref{eq:alspach}), we have:

\begin{corollary}\label{cor:main}
For $n\ge 4$,
\[{\rm cms}(K_n)=
\left\{\begin{array}{cl}
{\rm ms}(K_n)&\mbox{ if $n$ is even}\\
{\rm ms}(K_n)-1&\mbox{ if $n$ is odd.}\end{array}\right.\]  
\end{corollary}

\section{Proof of the Main Result}

The proof is in three parts: constructions for the even case (I), constructions for the   odd  case (II), and a nonexistence argument  for the odd case (III).

 \bigskip
 \begin{center}
\includegraphics{cms.30}
\qquad\includegraphics{cms.31}
\qquad\includegraphics{cms.32}

\bigskip
\includegraphics{cms.33}
\qquad\includegraphics{cms.34}

\end{center}
 
 \bigskip
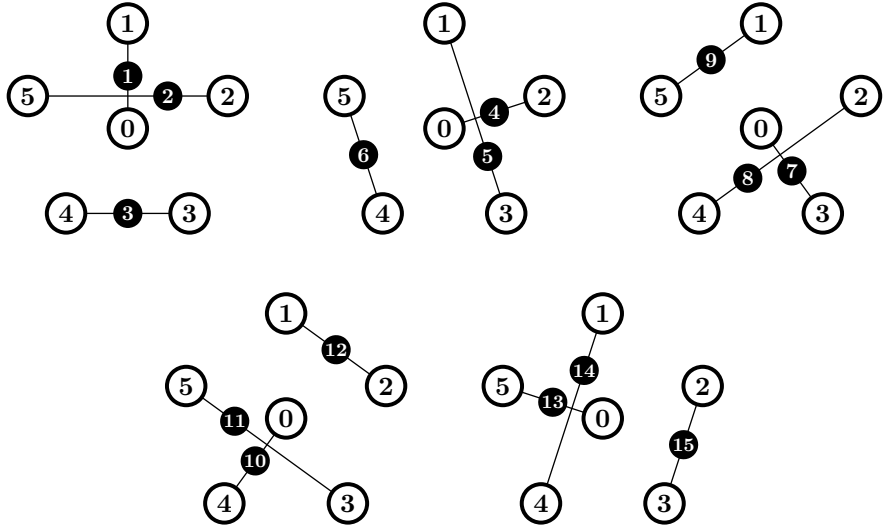
 \centerline{{\bf Figure 1}}
 \bigskip

(I) We first show  that if $m\ge 2$, then there is a cyclic  labeling of the edges of $K_{2m}$ whose cyclic matching number is $m-1$. In view of (\ref{eq:ineq}), this will establish  the first equality in Theorem \ref{th:main}.

Let the vertices of $K_{2m}$ be labeled $0,1,2,\ldots,2m-1$, 
and consider the vertices $1,2,\ldots,2m-1$ as equally spaced points on the unit circle centered at  a point labeled $0$. Consider the perfect matching 
\[M=\{\{0,1\},\{2,2m-1\}, \{3,2m-2\},\ldots, \{m,m+1\}\},\]
 of $K_{2m}$.
We label the edges of $M$ as listed in the order  $[1],[2],\ldots,[m]$.
 Let $\phi$ be the rotation 
about the origin by $360/(2m-1)$ degrees. Then 
\[ M, \phi(M),\phi^2(M),\ldots,\phi^{2m-2}(M)\]
is a 1-factorization of $K_{2m}$, that is,  $(2m-1)$ perfect matchings of $K_{2m}$ partitioning  the edges of $K_{2m}$. We note that $\phi^{2m-1}(M)=M$.
  For $k=1,2,\ldots,2m-2$, the edges
\[\phi^k(\{0,1\}),\phi^k(\{2,2m-1\}),\phi^k(\{3,2m-2\}),\ldots,\phi^k(\{m,m+1\})\]
receive labels 
\[[1+mk],[2+mk],[2+mk]\ldots,[m+mk],\]
 respectively.
 In this way we obtain an ordering $[1],[2],\ldots,[2m^2-m]$ of the edges of $K_{2m}$ which we regard as a cyclic ordering with [1] following $[2m^2-m]$. 
 To verify that this cyclic labeling has cyclic matching number  $m-1$,
 it suffices to show that each set of $m-1$ consecutive edges of $\phi^k(M),\phi^{k+1}(M)$  is a matching for each $k=1,2,\ldots,2m-2$.
 But this is straightforward to check using the  definitions of $M$ and $\phi$.
 We illustrate this construction  for $m=3$, that is, $K_6$, in Figure 1.

(II) We next show  that if $m\ge 2$, then there is a cyclic  labeling of the edges of $K_{2m+1}$ whose cyclic matching number is $m-1$.

A {\it near-perfect matching} of $K_{2m+1}$ is a matching of $m$ edges. Such a matching meets all but one vertex of $K_{2m+1}$, and we call the missing vertex the {\it isolated vertex of the matching}.
Let the vertices of $K_{2m+1}$ be labeled $0,1,2,\ldots,2m$, 
and consider the vertices $0,1,2,\ldots,2m$ as equally spaced points on the unit circle centered at the origin.
The near-perfect matching 
\[M=\{\{1,2m\}, \{2,2m-1\}, \ldots, \{m,m+1\}\}\]
has   vertex $0$ as its isolated vertex. We label the edges of $M$ as listed in the order $[1],[2],\ldots,[m]$ (according to the smallest vertex the edge contains). Let $\phi$ be the rotation 
about the origin by $360/(2m+1)$ degrees. Then 
\[ M, \phi(M),\phi^2(M),\ldots,\phi^{2m}(M)\]
are $(2m+1)$ near-perfect matchings of $K_{2m+1}$ and they partition the edges of $K_{2m+1}$. Each vertex of $K_{2m+1}$ is an isolated vertex of one of these near-perfect matchings.  For $k=1,2,\ldots,2m$, the edges
\[\phi^k(\{1,2m\}),\phi^k(\{2,2m-1\}),\ldots,\phi^k(\{m,m+1\})\]
receive labels 
\[[1+mk],[2+mk],\ldots,[m+mk],\]
 respectively.
 In this way we obtain an ordering $[1],[2],\ldots,[2m^2+m]$ of the edges of $K_{2m+1}$ which we regard as a cyclic ordering with [1] following $[2m^2+m]$.  (We remark that the Walecki decomposition of $K_{2m+1}$ into Hamilton cycles proceeds with $2m$ vertices equally spaced on the unit circle and one at the origin and  with one Hamilton cycle $C$. Then with $\theta$ equal to the rotation about the origin through $360/(2m)$ degrees $C,\theta(C),\ldots,\theta^{m-1}(C)$ is a Hamilton decomposition of $K_{2m+1}$.)

To verify that this cyclic labeling has cyclic matching number (at least) $m-1$, we first note that  each $\phi^k(M)$ is a near-perfect matching and thus  each of the two  sets of $m-1$ consecutive edges of $\phi^k(M)$ (that is, its first $m-1$ edges and its last $m-1$ edges) is a matching. Also since $\phi(\phi^k(M))= \phi^{k+1}(M)$ including $\phi(\phi^{2m} (M))=\phi^{2m+1}(M)=M$, it is enough to show that any $m-1$ consecutive edges from $M\cup \phi(M)$ form a matching.
But again this is straightforward to check using the  definitions of $M$ and $\phi$.  We illustrate this construction for $m=3$, that is, $K_7$, in Figure 2.

\bigskip
 
 \begin{center}
 \includegraphics{cms.20}
\qquad\includegraphics{cms.21}
\qquad\includegraphics{cms.22}

\bigskip
\qquad\includegraphics{cms.23}
\qquad\includegraphics{cms.24}
\qquad\includegraphics{cms.25}

\bigskip
\qquad\includegraphics{cms.26}

 \end{center}
 \bigskip
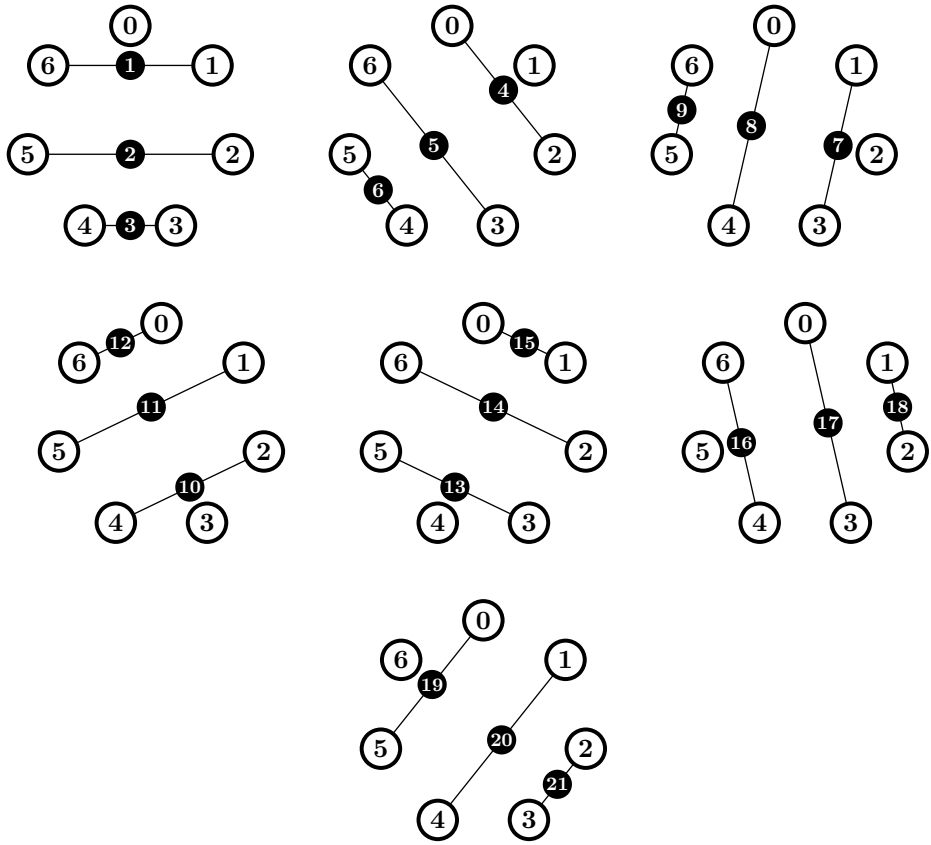
\centerline{{\bf  Figure 2}}
 \bigskip

(III) Finally we show that if $m\ge 2$, then there does not exist a cyclic ordering of the edges of $K_{2m+1}$  with cyclic matching number $m$.  In view of (II) we conclude that ${\rm cms}(K_{2m+1})=m-1$.

Assume to the contrary that there is a cyclic ordering $[1],[2],\ldots,[2m^2+m],[1]$ of the edges of $K_{2m+1}$ with cyclic matching number equal to $m$.  Each set of $m$ consecutive edges in the cyclic ordering is a near-perfect matching of $K_{2m+1}$ with a unique isolated vertex. Let
\[M_{[i]} =\{[i],[i+1],[i+2],\ldots,[i+m-1]\} \quad (i=1,2,\ldots,2m^2+m)
\]
where   the edges in the matching $M_{[i]}$ are taken in the cyclic ordering.  
Since
\[M_{[i+1]} =\{[i+1],[i+2],\ldots,[i+m-1], [i+m]\},\]
the edge $[i+m]$  joins the isolated vertex of $M_{[i]}$ to one of the vertices of edge $[i]$.
For each $i=1,2,\ldots,m$, we get a  partition 
\[{\mathcal F}_{[i]}= \{M_{[i]}, M_{[i+m]},M_{[i+2m]},\dots, M_{[i+(m-1)m]}\}\]
 of the edges of $K_{2m+1}$ into
$2m+1$ near-perfect matchings; here the subscripts are taken modulo
$2m^2+m$.  Since the degree of each vertex of $K_{2m+1}$ equals $2m$, each vertex of $K_{2m+1}$ is the isolated vertex of exactly one near-perfect matching of ${\mathcal F}_{[i]}$.

For each $i=1,2,\ldots,m$, let $G_{[i]}$ be the spanning subgraph of $K_{2m+1}$ consisting of  the $2m+1$ edges
\begin{equation}\label{eq:G_i}
\{[i],[i+m],[i+2m],\ldots,[i+2m^2]\}.\end{equation}
Note that $[i+2m^2+m]=[i]$ in our cyclic ordering, since $K_{2m+1}$ has exactly $2m^2+m$ edges. 
Each graph $G_{[i]}$ has exactly  $2m+1$ vertices and exactly $2m+1$ edges.

\smallskip\noindent
{\bf  Claim 1: }  Each $G_{[i]}$ is a connected unicyclic graph consisting of a cycle with pendent edges $($possibly none$)$ at  its vertices.

\smallskip
In fact, since the number of edges of $G_{[i]}$ equals the number of its vertices, if $G_{[i]}$ is connected and the edges span all its vertices,  $G_{[i]}$ must be  unicyclic. In constructing  $G_{[i]}$ in the order of the edges listed, each edge of $G_{[i]}$  has a vertex in common with the  edge of $G_{[i]}$ preceding 
it in the order (interpreted cyclically) given in (\ref{eq:G_i}); thus,
starting with one of the vertices of $G_{[i]}$,  we include one new vertex at a time, giving an ordering of the vertices for each $G_{[i]}$,   and attach an edge from the new vertex to an old vertex.  It follows that $G_{[i]}$ is connected, and since each vertex of $K_{2m+1}$ is an isolated vertex of  one of the near-perfect matchings of ${\mathcal F}_i$, the edges of $G_{[i]}$ span its vertices. Thus $G_{[i]}$ is a connected unicyclic graph, and so  $G_{[i]}$ consists of a cycle with a tree (possibly an empty tree) rooted at each of its vertices. Moreover, since each edge of $G_{[i]}$ has a vertex in common with the previous edge, these trees cannot have paths of length greater than 1. This establishes Claim 1.

Consider the connected unicyclic graph $G_{[1]}$. Let $k$ be the length of its cycle $\gamma$ where we  label the vertices so that $\gamma=(1,2,\ldots,k,1)$. We start the cyclic ordering of the edges of $G_{[1]}$ with the edge $\{1,2\}$, and we can assume that $\{1,2\}$ is edge [1] in the cyclic ordering of the edges of $K_{2m+1}$. Let there be $p_i\ge 0$ pendent edges of $G_{[1]}$ at vertex $i$.  Since each edge of $G_{[1]}$ has a vertex in common with the previous edge, it follows that the cyclic ordering of the edges of $G_{[1]}$ is of the form
\[\{1,2\}, \{2,v_{2,1}\},\ldots,\{2,v_{2,p_2}\}, \{2,3\},\{3,v_{3,1}\},\ldots,\{3,v_{3,p_3}\},
\ldots\] 
\[\ldots, \{k,1\},\{1,v_{1,1}\},\ldots,\{1,v_{1,p_1}\}, \{1,2\},\]
where $v_{i,1},\ldots,v_{i,p_i}$ are the pendent vertices of $G_{[i]}$ joined to vertex $i$ ($i=1,2,\ldots,k$).

Assuming, without loss of generality, that vertex 1 (as opposed to vertex $2$) is the isolated vertex of the near-perfect matching $M_{[2]}$, then it follows that the isolated vertices of the near-perfect matchings in the partition ${\cal F}_{[2]}$ are in the order 
\begin{equation}\label{eq:order}
1, v_{2,1},\ldots,v_{2,p_2}, 2,v_{3,1},\ldots,v_{3,p_3}
\ldots,k,v_{1,1},\ldots,v_{1,p_1}, 1.\end{equation}
Thus, in constructing $G_{[2]}$ in the  manner described above, the vertices are added in this order with an edge from a vertex to some preceding vertex.

\smallskip\noindent
{\bf  Claim 2:}  The edges $\{v_{t,i},v_{t,{i+1}}\}$ $(1\le i<p_t)$ and the edges $\{t-1,v_{t,1}\}$ $(1\le t\le k)$ are not edges of $G_{[2]}$. 

\smallskip
Suppose  $\{v_{t,i},v_{t,i+1}\}$ is an edge of $G_{[2]}$. Then $v_{t,i+2}$ 
must be joined to either $v_{t,i}$ or $v_{t,i+1}$ in $G_{[2]}$, and it follows that $v_{t,p_t}$ is joined by an edge to one of $v_{t,i},\ldots,v_{t,p_{t-1}}$.  Now vertex $t$ is joined by an edge in $G_{[2]}$ to one of
 $v_{t,i},\ldots,v_{t,p_{t-1}},v_{t,p_t}$, a contradiction since all such edges belong to $G_{[1]}$. 
 Similar reasoning shows that   $\{t-1,v_{t,1}\}$ cannot be  an edge of $G_{[2]}$. This establishes Claim 2.

The vertices of $G_{[2]}$ are added in the order given in (\ref{eq:order}) with its edge $[2+m]$ equal to $\{1,w\}$ for some vertex $w$. The edge $[2+2m]$ is the edge $\{v_{2,1},w\}$, since by Claim 2, it cannot be the edge $\{1,v_{2,1}\}$.  Similarly, by Claim 2, the next edge is  $\{v_{2,2},w\}$ since it cannot be the edge $\{v_{2,1},v_{2,2}\}$. Continued application of Claim 2 shows that all subsequent vertices must be joined by an edge to $w$ in $G_{[2]}$. This implies that $G_{[2]}$ is a graph $K_{1,2m}$ and thus is not unicyclic, a contradiction.

This completes the proof of Theorem \ref{th:main}.

\bigskip
\begin{center}
\includegraphics{cms.11}
\qquad\includegraphics{cms.14}

\bigskip
\includegraphics{cms.12}
\qquad\includegraphics{cms.15}

\bigskip
\includegraphics{cms.13}
\qquad\includegraphics{cms.16}
 
 \end{center}
 \bigskip
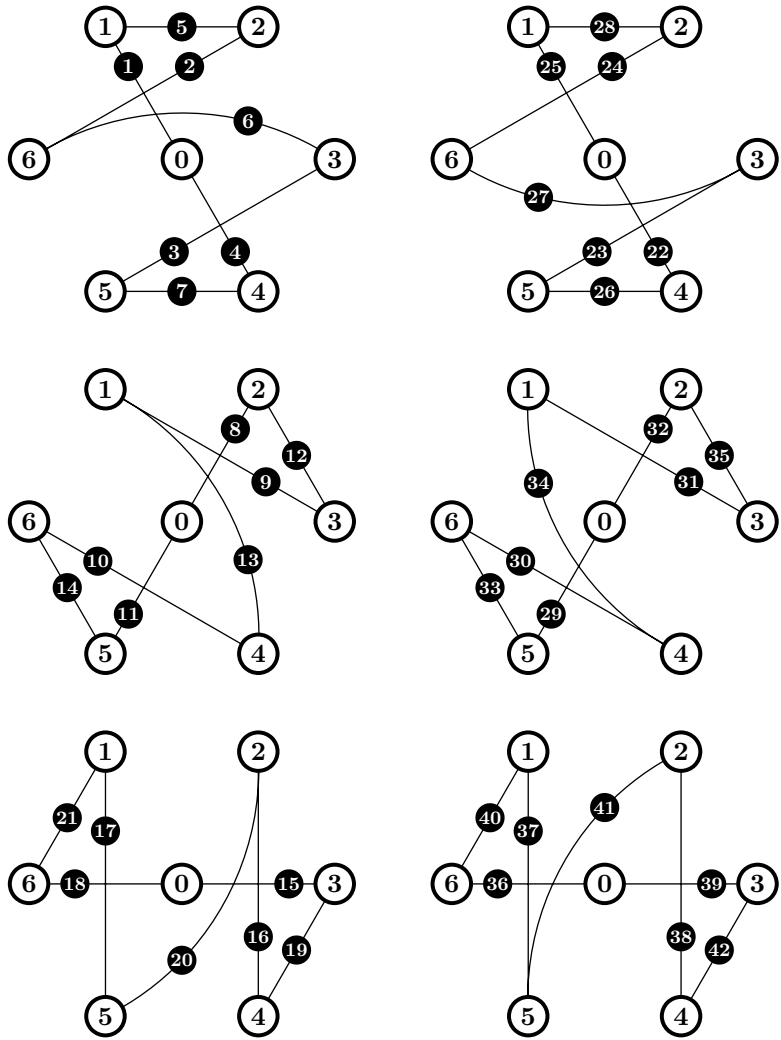
 \centerline{{\bf Figure 3}}
 \bigskip

\section{Remarks  and Open Questions}

Let $m\ge 2$.
By Theorem \ref{th:main}, ${\rm cms}(K_{2m+1})=m-1$. As already remarked, in \cite{Al}, the Walecki decomposition of $K_{2m+1}$ into Hamilton cycles is used to show  that ${\rm ms}(K_{2m+1})=m$.  If one uses that construction twice, one obtains that
\[{\rm cms}(2K_{2m+1})=m\quad (m\ge 2),\]
where $2K_{2m+1}$ is the multigraph obtained from $K_{2m+1}$
by doubling each edge. This is illustrated in Figure 3 for $m=3$. In this figure, the ordering of the edges of $K_7$ given by the 
integers $1,2,\ldots, 21$ is the ordering from the Walecki decomposition used by Alspach  \cite{Al} to show that the matching sequencibility  of $K_7$ 
equals 3; note that this ordering does not give that the cyclic matching sequencibility of $K_7$ equals 3.

If $C_n$ is a cycle of $n$ vertices with $n$ odd, then it is easy to show that ${\rm cms}(C_n)={\rm ms}(C_n)=\frac{n-1}{2}$, the maximum possible. This follows by assigning the  integers $1,2,\ldots,n$ to the edges, starting with any edge and going around the cycle twice taking alternate edges. This is illustrated in Figure 4 for $n=7$.

\bigskip
\begin{center}
\includegraphics{cms.41}

 \end{center}

 \bigskip
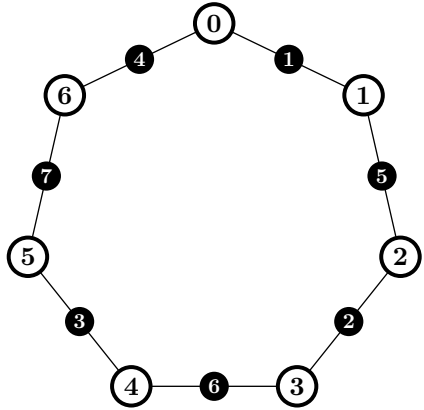
 \centerline{{\bf Figure 4}}
 
 \bigskip
 
 Now consider an even length cycle $C_{2q}$. Then $C_{2q}\subseteq K_{q,q}$ and $C_{2q}$
has as a biadjacency matrix $A=I+P_q$ where $P_q$ is the permutation
matrix with 1s in positions $(1,2), (2,3),\ldots, (q-1,q), (q,1)$.

First suppose that $q$ is even, and let a cyclic ordering of the edges of $C_{2q}$ be given by the matrix $\widehat{A}$ equal to
\[
\left[\begin{array}{c|c|c|c|c|c|c|c|c|c|c|c|c}
1&q&&&&&&&&&&&\\ \hline
&2q-1&q-2&&&&&&&&&&\\ \hline
&&2q-3&q-4&&&&&&&&&\\ \hline
&&&2q-5&q-6&&&&&&&&\\ \hline
&&&&\ddots&\ddots&&&&&&&\\ \hline
&&&&&q+3&2&&&&&&\\ \hline
&&&&&&q+1&2q&&&&&\\ \hline
&&&&&&&q-1&2q-2&&&&\\ \hline
&&&&&&&&\ddots&\ddots&&&\\ \hline
&&&&&&&&&9&q+8&&\\ \hline
&&&&&&&&&&7&q+6&\\ \hline
&&&&&&&&&&&5&q+4\\ \hline
q+2&&&&&&&&&&&&3\end{array}\right].
\]

Now the difference of two positive integers in the same row or column of
$\widehat{A}$ is $\pm (q-1)$ modulo $2q$.  This implies that every (cyclical) set of
$q-1$ consecutive positive  integers lie in different rows and columns.  Since ${\rm
cms}(C_{2q})$ and ${\rm ms}(C_{2q})$ are at most $q-1$, it follows that \[{\rm
cms}(C_{2q})={\rm ms}(C_{2q})=q-1 \quad (q\mbox{ equal $0$ modulo $4$}).\] The matrix $\widehat{A}$ when $q=4$ is 
 \[
\left[\begin{array}{c|c|c|c|c|c|c|c} 1&8&&&&&&\\ \hline &15&6&&&&&\\ \hline
&&13&4&&&&\\ \hline &&&11&2&&&\\ \hline &&&&9&16&&\\ \hline &&&&&7&14&\\ \hline
&&&&&&5&12\\ \hline 10&&&&&&&3\end{array}\right].\]

Now suppose that $q$ is odd,  and let the cyclic ordering of the edges of $C_{2q}$  be given by the matrix 
$\widehat{A}$ equal
to
\[
\left[\begin{array}{c|c|c|c|c|c|c|c|c|c|c|c|c}
1&q+2&&&&&&&&&&&\\ \hline
&3&q+4&&&&&&&&&&\\ \hline
&&5&q+6&&&&&&&&&\\ \hline
&&&7&q+8&&&&&&&&\\ \hline
&&&&\ddots&\ddots&&&&&&&\\ \hline
&&&&&q-1&2q&&&&&&\\ \hline
&&&&&&q&2q-1&&&&&\\ \hline
&&&&&&&q-2&2q-3&&&&\\ \hline
&&&&&&&&\ddots&\ddots&&&\\ \hline
&&&&&&&&&8&q+7&&\\ \hline
&&&&&&&&&&6&q+5&\\ \hline
&&&&&&&&&&&4&q+3\\ \hline
q+1&&&&&&&&&&&&2\end{array}\right].
\]
The difference of two positive integers in the same row and column of $\widehat{A}$ is $\pm (q-1)$ or $q$ modulo $2q$ implying that  every (cyclical) set of $q-1$ consecutive positive integers  lies in different rows and columns. Thus 
\[{\rm.
cms}(C_{2q})={\rm ms}(C_{2q})=q-1 \quad (q\mbox{ equal 2 modulo $4$}).\]
The matrix $\widehat{A}$ when $q=6$ is 
\[\left[\begin{array}{c|c|c|c|c|c}
1&8&&&&\\ \hline
&3&10&&&\\ \hline
&&5&12&&\\ \hline
&&&6&11&\\ \hline
&&&&4&9\\ \hline
7&&&&&2\end{array}\right].\]
\bigskip

In summary, we have
\[{\rm cms}(C_n)={\rm ms}( C_n)=\left\lfloor\frac{n-1}{2}\right\rfloor\quad (n\ge 3).\]

Now consider a path $P_n$ of $n$ vertices.  First assume that  $n=2q$.  
Then a $q\times q$  biadjacency matrix $A$ of $P_{2q}$ is obtained from the biadjacency matrix of the cycle $C_{2q}$ by replacing the 1 in its $(q,1)$-position with a zero.
It is easy to obtain an ordering of the edges of $P_{2q}$ (the $1$s in $A$)  that shows that ${\rm cms}(P_{2q})={\rm ms}(P_{2q})=q-1$. For example, if $q=5$, the following ordering works and generalizes in the obvious way:
\[\widehat{A}=\left[\begin{array}{c|c|c|c|c}
4&8&&&\\ \hline
&3&7&&\\ \hline
&&2&6&\\ \hline
&&&1&5\\ \hline
&&&&9\end{array}\right].\]

Now assume that $n=2q+1$ is odd. 
Then a $q\times (q+1)$  biadjacency matrix $A$ of $P_{2q+1}$ is obtained from the biadjacency matrix of the cycle $C_{2q+1}$ by deleting its last row. It is easy to show that ${\rm ms}(P_{2q+1})=q$. For example, if  $q=5$,  the following ordering works and generalizes in the obvious way:
\[\widehat{A}=\left[\begin{array}{c|c|c|c|c|c}
5&10&&&&\\ \hline
&4&9&&&\\ \hline
&&3&8&&\\ \hline
&&&2&7&\\ \hline
&&&&1&6\end{array}\right].\]
Note that this construction does not give ${\rm cms}(P_{2q+1})=q$. In fact, ${\rm cms}( P_{2q+1})
=q-1$, as we now argue by contradiction. Suppose there were a cyclic ordering of the $2q$ edges of $P_{2q+1}$ such that every cyclic set of $q$ edges form a matching of $P_{2q+1}$.  Since we have a cyclic ordering, we may assume that the edge $[1]$ is one of the pendent edges of $P_{2q+1}$. But then, no matter which matching is determined by the edges $[1],[2],\ldots,[q]$, all of the other edges
meet a vertex of one of the edges $[2],\ldots,[q]$, a contradiction.

In summary, we have
\[\begin{array}{ccl}
{\rm cms}(P_n)={\rm ms}(P_n)&=&\frac{n-2}{2}\quad \mbox{ if $n$ is even,}\\
{\rm ms}(P_n)&=&\frac{n-1}{2}\quad \mbox{ if $n$ is odd,}\\
{\rm cms}(P_n)&=&\frac{n-3}{2}\quad \mbox{ if $n$ is odd.}\end{array}\]

The (cyclic) matching sequencibility has been computed for complete graphs, complete bipartite graphs, cycles, and paths. It may be of interest
to compute it for other important  classes of graphs, but such computations may be very difficult and they will in general depend on the particular graph in the class and thus be of less interest. 
For a tree,  the (cyclic) matching sequencibility may differ considerably from the  {\it matching number} (the maximum number of edges in a matching). For example, in Figure 5 there is a tree $T$ 
of order $18$ with a perfect matching (matching number equal to 9)  and an ordering of the edges showing that ${\rm ms }(G)\ge 2$.  It is not hard to see that ${\rm cms}(T) =1$ and ${\rm ms}(T)=2$.  If one takes any graph $G$ of order   $n$  and adjoins $n+1$ (respectively, $n+2$) pendent edges at one of its  vertices to get a graph $H$ of order $2n+1$ (respectively, $2n+2$), then ${\rm ms}(G) =1$ (respectively,  ${\rm cms}(H)=1$), since in any ordering  (respectively, cyclic ordering) of the edges of $H$,  some two of these  pendent edges would have to be consecutive.

\bigskip
\begin{center}
\includegraphics{cms.51}
\end{center}

\bigskip
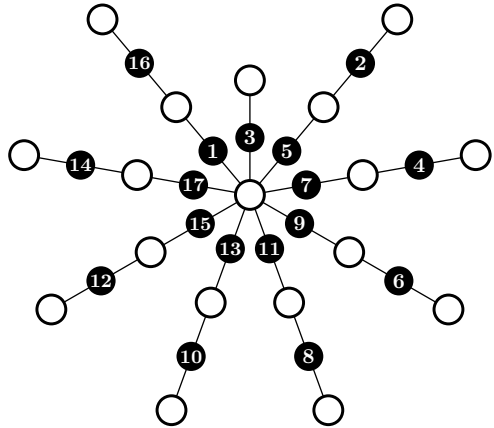
\centerline{\bf Figure 5}

\bigskip

Another class of graphs of potential interest are the $k$-regular bipartite graphs $G\subseteq K_{n,n}$.
For example, consider the 3-regular bipartite graph $G\subseteq K_{n,n}$ whose biadjacency matrix is $P^{-1}+I +P$ 
(equivalently, $I+P+P^2$) where $P$ is as before the full cycle permutation matrix of order $n$. Then ${\rm cms}(G)={\rm ms }(G)=n-1$. We illustrate this in terms of the biadjacency matrix $A$ and the labeling given by $\widehat{A}$ for the $n$ odd  and $n$ even cases with $n=7$ and $8$:
\[
\widehat{A}=
\left[\begin{array}{c|c|c|c|c|c|c}
1&16&&&&&7\\ \hline
15&2&9&&&&\\ \hline
&10&3&18&&&\\ \hline
&&17&4&11&&\\ \hline
&&&12&5&20&\\ \hline
&&&&19&6&13\\ \hline
8&&&&&14&21\end{array}\right],\quad
\left[\begin{array}{c|c|c|c|c|c|c|c}
1&18&&&&&&8\\ \hline
17&2&10&&&&&\\ \hline
&11&3&20&&&&\\ \hline
&&19&4&12&&&\\ \hline
&&&13&5&22&&\\ \hline
&&&&21&6&14&\\ \hline
&&&&&15&7&24\\ \hline
9&&&&&&23&16\end{array}\right].\]

We conclude with three additional questions applying to all graphs that appear to be very difficult.

\smallskip

\noindent
{\bf Question 1:}  Given a graph $G$ with matching number  $p$, is there a positive integer $k$ such that 
${\rm ms}(kG)=p$
(${\rm cms}(kG)=p$)?

\medskip\noindent
{\bf Question 2:}  For a graph $G$, we have ${\rm ms}(G)\ge {\rm cms}(G)$. How large can ${\rm ms}(G)-{\rm cms}(G)$ be? Is ${\rm cms}(G)\ge {\rm ms}(G)-1$?  

\medskip\noindent
{\bf Question 3:} Given a graph $G$, is ${\rm cms}(2G)={\rm ms}(G)$?

\bigskip\noindent


\begin{thebibliography}{99}
\bibitem{Al}  B.~Alspach, The wonderful Walecki construction,  {\it Bull. of the ICA}, 52 (2008), 7--20.



\end{thebibliography}
\end{document}